\documentclass{article}
\usepackage{graphicx} 

\usepackage{amsmath}
\usepackage{amssymb}
\usepackage{amsfonts}
\usepackage{amsthm}
\usepackage{mathtools}
\usepackage{enumerate}
\usepackage[hidelinks]{hyperref}
\usepackage{titlesec}
\usepackage{url}

\usepackage[T1]{fontenc}   
\usepackage{lmodern}       
\usepackage{textcomp}      

\usepackage[shortlabels]{enumitem}
\usepackage[english]{babel}
\usepackage{yfonts}
\usepackage{setspace}
\usepackage[utf8]{inputenc}
\usepackage[english]{babel}
\usepackage[margin = 1in]{geometry}
\usepackage{tikz}
\usetikzlibrary{cd}
\usepackage{graphicx}
\usepackage{authblk}
\usetikzlibrary{calc}
\usepackage{fancyhdr}
\usepackage{float}

\theoremstyle{definition}

\theoremstyle{remark}
\newtheorem*{rem}{Remark}
\theoremstyle{plain}
\newtheorem{conj}{Conjecture}

\newtheorem{thm}{Theorem}[section]
\newtheorem{prop}[thm]{Proposition}
\newtheorem{lem}[thm]{Lemma}
\newtheorem{cor}[thm]{Corollary}
\newtheorem{defi}[thm]{Definition}

\titleformat{\section}[block]{\centering\Large\bfseries}{\thesection}{1em}{}

\title{The Table Problem for Convex Domains}
\author{Ali Naseri Sadr}
\date{}

\begin{document}

\maketitle
\begin{abstract}
    We prove a generalized Table Theorem for bounded convex domains in the plane. Then we prove a variation of the Table Problem for compactly supported functions on the plane. Our method is based on relating the table problem on a surface to the square peg problem for a configuration of star-shaped curves in the tangent bundle of the surface.
\end{abstract}

\section{Introduction}

The Table Problem suggested by Knaster in $1947$ conjectures: 
\begin{conj}[The Table Problem, Problem $2.14$ in \cite{MR5077826}]
\label{table conj}
    For every continuous function $f\colon S^2\to \mathbb{R}$ and $d\in(0,\sqrt{2}]$, there are four points $x_1, x_2, x_3, x_4$ on the unit sphere in $\mathbb{R}^3$ such that they are vertices of a square with side length $d$ and
    \begin{equation*}
        f(x_1) = f(x_2) = f(x_3) = f(x_4).
    \end{equation*}
\end{conj}
The idea behind this problem is to think of the function $f$ as a height function on the two sphere, and the points $x_1,\cdots, x_4$ as the basis of a square shaped table of size $d$. In other words, the problem asks if we can always balance a square shaped table on a planet. 

Dyson proved the case for $d=\sqrt{2}$ in \cite{Dyson}, and Fenn proved a variation of this problem for non-negative functions that are supported inside a bounded convex domain in \cite{Fenn}. We generalize Fenn's result to every function that is constant outside of a bounded convex domain in the plane. The non-negativity is essential to Fenn's argument. Our argument is different from Fenn, and we are able to remove the non-negativity condition by a method based on our proof of Conjecture \ref{table conj} in \cite{sadr2026tabletheoremsurfacesodd} for even functions. 

\begin{defi}
    Consider a continuous function $f\colon\mathbb{C}\to \mathbb{R}$. A table of radius $a$ for $f$ is a square of radius $a$ in $\mathbb{C}$ such that $f$ takes the same value over the four vertices of the square.
\end{defi}
\begin{thm}[\textbf{Table Theorem for Convex domains}]
\label{Table for convex domains}
    Let $\Omega$ be a bounded convex domain in $\mathbb{C}$, and consider a continuous function $h\colon\mathbb{C}\to \mathbb{R}$ that is constant outside of $\Omega$. Then $h$ admits a square table of radius $a$ centered in $\Bar{\Omega}$ for every positive $a$.
\end{thm}

\begin{rem}
    This result is false if we remove the convexity assumption on $\Omega$; see \cite{MR638185} for a counterexample.
\end{rem}
Theorem \ref{Table for convex domains} leads to the following variation of Conjecture \ref{table conj} for compactly supported continuous functions on the plane. 
\begin{cor}[\textbf{The Planar Table Theorem}]
\label{Planar Table Thm}
  Let $f$ be a compactly supported continuous function on the plane. Then $f$ admits a table of radius $a$ centered in the convex hull of $\operatorname{supp}(f)$ for every $a>0$.  
\end{cor}
\begin{proof}
    Take $\Omega$ to be the convex hull of $\operatorname{supp}(f)$ and apply Theorem \ref{Table for convex domains} to $\Omega$ and $f$.
\end{proof}

\section*{Acknowledgments}
 The author is grateful to his advisors, John Baldwin and Josh Greene, for their invaluable guidance, support, and insightful conversations about this work.

\section{Smooth Domains and Configuration of Star-shaped Curves}

Suppose $\Omega$ is a bounded convex domain in $\mathbb{C}$ with $C^1$ boundary. We denote the space of positive $C^3$ functions that are constant outside of $\Omega$ by $\mathcal{F}(\Omega)$. Fix a positive number $a$, and let $N_{2a}(\Omega)$ be the open neighborhood of $\Omega$ in $\mathbb{C}$ that consists of all the points with distance less than $2a$ from $\Omega$. Define $\mathcal{F}(N_{2a}(\Omega))$ in a similar manner. Both $\mathcal{F}(\Omega)$ and $\mathcal{F}(N_{2a}(\Omega))$ are Banach manifolds with the $C^3$ norm. Let $X$ be a topological space and denote its fourth power by $X^4$. We can permute the coordinates in $X^4$, and this induces an action of $S_4$ on this space. 
\begin{defi}
 Define the fourth symmetric product of $X$ by
    \begin{equation*}
        \operatorname{sym}^{4}(X)\coloneqq X^4/S_4.
    \end{equation*}
\end{defi}

We define the fat diagonal $\Delta$ to be the subset of points in $\operatorname{sym}^4(X)$ for which at least two of the coordinates are equal. Define $K_a(\Bar{\Omega})$ to be a fiber bundle over $\Bar{\Omega}$ where each fiber over a point $z$ is the fourth symmetric product of the circle with radius $a$ in $T_z\mathbb{C}$ minus its fat diagonal; observe that $K_a(\Bar{\Omega})$ is a non-orientable eight dimensional manifold with boundary. Denote $\operatorname{sym}^4(\mathbb{C})\setminus\Delta$ by $V$ and let $A$ be the subset of four tuples in $V$ that are vertices of a square; this is a four dimensional submanifold of $V$. We define a map $\Psi\colon \mathcal{F}(N_{2a}(\Omega))\times K_a(\Bar{\Omega})\to \Bar{\Omega}\times V$ by
\begin{equation*}
    (h, (x, [v_1, v_2, v_3, v_4]))\mapsto (x, [h(x+v_1)\cdot v_1, h(x+v_2)\cdot v_2, h(x+v_3)\cdot v_3, h(x+v_4)\cdot v_4]).
\end{equation*}
Note that if we fix a function $h\in \mathcal{F}(N_{2a}(\Omega))$, then each fiber of $K_a(\Bar{\Omega})$ is mapped to the fourth symmetric product of a star-shaped curve in the tangent bundle of $\Bar{\Omega}$ by $\Psi$; we denote this map by $\Psi_h\colon K_a(\Bar{\Omega})\to \Bar{\Omega}\times V$, and we call $\Psi_h(K_a(\Bar{\Omega}))$ the configuration of star-shaped curves associated to $h$.
\begin{lem}
    The map $\Psi$ is a submersion.
\end{lem}
\begin{proof}
     Choose a pair $(h, (x,[v_1, v_2, v_3, v_4]))$ in $\mathcal{F}(N_{2a}(\Omega))\times K_a(\Bar{\Omega})$. Observe that since $\Psi_h$ maps the fiber over $x$ in $K_a(\Bar{\Omega})$ to the fiber over $x$ in $\Bar{\Omega}\times V$, the map $d\Psi$ is an isomorphism over the horizontal vectors, and it suffices to show the vertical vectors are in the image of $d\Psi$. Fix all the $v_i's$ except $v_1$, and change $v_1$ along an arc $\delta(t)$ in the circle with radius $a$ around the origin in $T_x\mathbb{C}$ such that the curve $\gamma(t) = [\delta(t),\theta_2,\theta_3,\theta_4]$ avoids the fat diagonal in $sym^4(\mathbb{C})$. We have
    \begin{equation*}
        \Psi(h,(x,\gamma(t))) = (x,[h(x+\delta(t))\delta(t),h(x+v_2)v_2,h(x+v_3)v_3,h(x+v_4)v_4]).
    \end{equation*}
    Hence, we get 
    \begin{equation*}
        d\Psi(0,(0,\Dot{\gamma}(0))) = [h(x+v_1)\cdot\Dot{\delta}(0) +\frac{d(h(x+\delta(t)))}{dt}\Big\vert_{t=0}\cdot v_1,0,0,0] \in \operatorname{Im}(d\Psi),
    \end{equation*}
    where we used $\delta(0) = v_1$. Since $\delta(t)$ has constant length, $\Dot{\delta}(0)$ is a non-zero vector orthogonal to $v_1$. Repeating this argument for the other coordinates proves we have vectors of the previous form in the image of $d\Psi$ where all the components are zero except one and the non-zero component is of the form $r_i\cdot \omega_i+s_i\cdot v_i$ with $r_i = h(x+v_i)$ greater than zero, $\omega_i$ a vector orthogonal to $v_i$, and $s_i$ some arbitrary real number. Now fix the four-tuple $v$ and pick a bump function $g$ supported in $N_{2a}(\Omega)$ such that $g(x+v_1) = 1$ and $g$ vanishes on the other $x+v_i$'s. Let $g_t = h+tg$ and note that for sufficiently small $t$, all the functions $g_t$ are in $\mathcal{F}(N_{2a}(\Omega))$. We get
    \begin{equation*}
        d\Psi(g,(0,0)) = [g(x+v_1)\cdot v_1,g(x+v_2)\cdot v_2,g(x+v_3)\cdot v_3,g(x+v_4)\cdot v_4] = [v_1,0,0,0]\in \operatorname{Im}(d\Psi) .
    \end{equation*}
 The same argument shows we have vectors of the form $[0,\dots,v_i,\dots,0]$ in the image of $d\Psi$. We conclude the lemma because the set of vectors
 \begin{align*}
     &[0,\dots,v_i,\dots,0], \\ 
     & [0,\dots,\omega_i,\dots,0]
 \end{align*}
 generate the vertical vectors over $x$ in $\Bar{\Omega}\times V$.
\end{proof}

\begin{rem}
    If we restrict $\Psi$ to $\mathcal{F}(\Omega)\times K_a(\Bar{\Omega})$, then it is not a submersion because we cannot move $h$ in the direction of the bump function unless $x+v_i$ is inside $\Omega$.
\end{rem}
Let 
\begin{equation*}
    \mathcal{S}\coloneqq \Psi^{-1}(\Bar{\Omega}\times A).
\end{equation*}
We call $\mathcal{S}$ the space of star-shaped squares over $\Bar{\Omega}$. This is a $C^3$ Banach manifold with boundary, and its codimension is $4$. Denote the first projection map on $\mathcal{F}(N_{2a}(\Omega))\times K_a(\Bar{\Omega})$ by $pr_1$, and let $\pi$ be its restriction to $\mathcal{S}$.

\begin{lem}
The map $\pi\colon \mathcal{S}\to \mathcal{F}(N_{2a}(\Bar{\Omega}))$ is Fredholm, and its index is $2$.    
\end{lem}
\begin{proof}
    Let $\iota\colon \mathcal{S}\to \mathcal{F}(N_{2a}(\Bar{\Omega}))\times K_a(\Bar{\Omega})$ denote the inclusion map. Both $\iota$ and $pr_1$ are Fredholm, and we have
    \begin{equation*}
        \pi = pr_1\circ\iota
    \end{equation*}
    so $\pi$ is Fredholm as well, and we get
    \begin{equation*}
        \operatorname{ind}(\pi) = \operatorname{ind}(pr_1)+\operatorname{ind}(\iota) = 6-4=2.
    \end{equation*}
\end{proof}

The map $\pi$ is of class $C^3$, and it has index $2$ so its regular values are dense in $\mathcal{F}(N_{2a}(\Omega))$ by Sard-Smale theorem. Denote the regular values of $\pi$ by $\mathcal{R}_{\Omega}$. Note that a function $h$ is in $\mathcal{R}_{\Omega}$ iff the intersection
\begin{equation*}
    \Psi_h(K_a(\Bar{\Omega}))\cap (\Bar{\Omega}\times A)
\end{equation*}
is transverse. This subspace corresponds to the set of inscribed squares in the configuration of star-shaped curves associated to $h$. In particular, we know that it is non-empty and compact for every $h\in \mathcal{F}(N_{2a}(\Omega))$ by the square peg problem for $C^2$ curves; we refer the unfamiliar reader to \cite{Matschke2011EquivariantTM, SHN,MR3184501}. 

\begin{defi}
    We define the center map $c\colon\mathcal{S}\to \mathbb{C}$ by 
    \begin{equation*}
        c\big(h, (x, [v_1,v_2,v_3,v_4])\big)\mapsto \frac{\sum_i h(x+v_i)\cdot v_i}{4}. 
    \end{equation*}
\end{defi}
Observe that if we have a point $(h, (x,[v_1, v_2,v_3, v_4]))\in \mathcal{S}$ with vanishing center, then we must have 
\begin{equation*}
    h(x+v_1) = h(x+v_2) = h(x+v_3) = h(x+v_4),
\end{equation*}
and the points $(x+v_1, x+v_2, x+v_3, x+v_4)$ are vertices of a square of radius $a$ centered at $x$. In other words, these points are the vertices of a table with radius $a$ for $h$.

\begin{prop}
\label{Tranvserse or vanishing on the boundary}
    Suppose $\Omega$ is a bounded convex domain in $\mathbb{C}$ with $C^1$ boundary. Then for every function $h\in \mathcal{F}(\Omega)$ and positive number $a$, one of the followings holds.
    \begin{enumerate}
        \item There is table of radius $a$ for $h$ with center on $\partial\Omega$.
        \item The center map restricted to 
        \begin{equation*}
            \Psi_h(K_a(\partial \Omega))\cap (\partial\Omega\times A)
        \end{equation*}
        maps each point to a vector that is transverse to $\partial\Omega$.
    \end{enumerate}
\end{prop}
\begin{proof}
    Assume that there is no table for $h$ with radius $a$ and center on $\partial\Omega$. Then $c$ restricted to 
    \begin{equation*}
        \Psi_h(K_a(\partial \Omega))\cap (\partial\Omega\times A)
    \end{equation*}
    must be non-vanishing. Choose a point $(h, (x,v_1, v_2, v_3, v_4))$ in this intersection, and denote the vector tangent to $\partial\Omega$ at $x$ by $w_x$ and the outward pointing normal by $n_x$. By contradiction, suppose $p\coloneqq c(h, (x,v_1, v_2, v_3, v_4))$ lies on the line generated by $w_x$, and it is not zero. We know the points $(h(x+v_i)\cdot v_i)$ are vertices of a square, and they lie on the star-shaped curve defined by
    \begin{equation*}
        v\mapsto h(x+v)\cdot v.
    \end{equation*}
    Since $\Omega$ is convex, $x+v$ is outside of $\Omega$ for every $v$ with $\langle v,n_x\rangle>0$. In particular, for every $v$ with $\langle v,n_x\rangle>0$ $h(x+v)$ is a positive constant by assumption. Therefore, half of our star-shaped curve is a half-circle centered around the origin, and the other half varies according to the values of $h$ inside $\Omega$. Denote each $h(x+v_i)\cdot v_i$ by $p_i$; these are vertices of a square, so at most three of them can lie on the half circle part of the star-shaped curve. If three of them are in the half circle part, then two of them are opposites of each other, and that means the center vanishes. Hence, assume at most two of them are in the half-circle part. If there is only one in the half-circle part, then we can find two opposite vertices $p_1$ and $p_3$ such that 
    \begin{align*}
       & \langle p_1, n_x\rangle <0, \\
       & \langle p_3, n_x\rangle<0.
    \end{align*}
    The center is the midpoint of $p_1$ and $p_3$, so we must have $\langle c, n_x\rangle<0$ which is a contradiction. The only remaining case is when two of the vertices are in the half-circle part, and the other two are not. Take the two vertices that are on the half-circle part; we know the center lies on the bisector of the line segment between these vertices. Denote this bisector by $l$; see the following figure. Moreover, we know that the center lies on the line generated by $w_x$ by assumption. Note that the diagonal of our half-circle is also the line generated by $w_x$, and the only intersection point between this line and $l$ is the origin. Thus the center must be the origin which gives a contradiction.  
    \begin{figure}[h]
        \centering
        \begin{tikzpicture}[scale=1]

    \draw[red, thick, ->] (-3,0) -- (3,0) node[right] {$w_x$};

    \draw[thick, blue]
        (2,0) arc[start angle=0, end angle=180, radius=2];

    \draw[thick, blue, domain=180:360, samples=200, variable=\t]
        plot (
            {(2 + 0.7*(sin(\t))^6*(1+0.6*sin(3*\t)))*cos(\t)},
            {(2 + 0.7*(sin(\t))^6*(1+0.6*sin(3*\t)))*sin(\t)}
        );

    \coordinate (O) at (0,0);
    \fill (0,0) circle (1.2pt);
    \node[below right] at (0,0) {$O$};

    \coordinate (A) at ({2*cos(110)},{2*sin(110)});
    \coordinate (B) at ({2*cos(160)},{2*sin(160)});

    \fill (A) circle (1.2pt);
    \fill (B) circle (1.2pt);

    \draw[dashed] (A) -- (B);

    \coordinate (M) at ($ (A)!0.5!(B) $);
    \fill (M) circle (1pt);
    \node[above right] at (-0.5,0.5) {$l$};

    \draw[dashed] (M) -- (O);

    \node[blue] at (1.7,-1.6) {$\gamma$};

\end{tikzpicture}
    \caption{The star-shaped curve $\gamma$ in $T_x\mathbb{C}$}
    \label{star-shaped curve in the tangent}
\end{figure}

\end{proof}
\section{Proof of the Main Theorem}
\begin{proof}[\textbf{Proof of Theorem \ref{Table for convex domains}}]
    First assume $\Omega$ has $C^1$ boundary. It suffices to prove the theorem for positive $C^3$ functions because if the result holds for a function $f$, then it also holds for every $f+c$, so we can shift every function by a large constant to make it positive; note this is possible since we assume $f$ is constant outside of $\Omega$. Take a function $h\in \mathcal{F}(\Omega)$
    Fix a positive $a$; if $h$ admits a table of radius $a$ centered on $\partial\Omega$, we get our claim. Thus suppose $h$ does not have a table centered on a boundary point; then we know the center map restricted to
    \begin{equation*}
            \Psi_h(K_a(\partial \Omega))\cap (\partial\Omega\times A)
    \end{equation*}
    maps each point to a vector that is transverse to $\partial\Omega$ by Proposition \ref{Tranvserse or vanishing on the boundary}. We can find a sequence of functions $h_n$ in $\mathcal{R}_{\Omega}$ that converges to $h$ in $C^3$ because $\mathcal{R}_\Omega$ is dense in $\mathcal{F}(N_{2a}(\Omega))$. For each $h_n$, the intersection
    \begin{equation*}
        \Sigma_n\coloneqq \Psi_{h_n}(K_a(\Bar{\Omega}))\cap \Bar{\Omega}\times A\subset \Bar{\Omega}\times V
    \end{equation*}
    is an embedded surface with boundary. This surface is properly embedded in $\Bar{\Omega}\times V$, and the first projection map to $\Bar{\Omega}$ has non-vanishing mod $2$ degree by the square peg problem for $C^2$ curves. Therefore, each $[\Sigma_n]$ is a non-zero element of
    \begin{equation*}
        H_2(\Bar{\Omega}\times V, \partial\Omega\times V; \mathbb{Z}_2)\cong H_2(\Bar{\Omega}, \partial\Omega;\mathbb{Z}_2)
    \end{equation*}
    The center map $c\colon \Bar{\Omega}\times V\to \Bar{\Omega}\times \mathbb{C}$ is a homotopy equivalence, and $[c(\Sigma_n)]$ is a non-zero element of
    \begin{equation*}
        H_2(\Bar{\Omega}\times\mathbb{C}, \partial\Omega\times\mathbb{C};\mathbb{Z}_2)\cong H_2(T\Omega, \partial T\Omega;\mathbb{Z}_2).
    \end{equation*}
    For sufficiently large $n$, $c(\Sigma_n)$ is transverse to the boundary of $\Omega$ because $h_n\to h$. Hence, the mod $2$ intersection number of $[c(\Sigma_n)]$ and $[\Bar{\Omega}]\times \{0\}$ is
    \begin{equation*}
        e(TD^2, \partial TD^2) =  \chi(D^2) = 1 \hspace{3mm} \text{mod $2$}.
    \end{equation*}
    We conclude that there is a $p_n$ in the interior of each $\Sigma_n$ with $c(p_n)=0$, and a limit point of these points lies in
    \begin{equation*}
        \Psi_h(K_a(\Bar{\Omega}))\cap (\Bar{\Omega}\times A).
    \end{equation*}
    This limit point gives a table of radius $a$ for $h$. Now suppose $\Omega$ is a general convex domain, and $h$ is a continuous function that is constant outside of $\Omega$. Approximate $\Omega$ by a sequence of convex domains $\Omega_n$ that contain $\Omega$ and have $C^1$ boundary. We know $h$ admits a table of radius $a$ with center in $\Bar{\Omega}_n$ for every $a$. A limit point of these tables for each $a$ is the desired table for $h$ with center in $\Bar{\Omega}$ and radius $a$. 
\end{proof}
Our proof of Theorem \ref{Table for convex domains} is based on relating the square peg problem to the table problem; inspired by the rectangular peg problem for smooth curves (we refer the unfamiliar reader to \cite{RectanglePaper}), we conjecture the following generalization of Theorem \ref{Table for convex domains}. 
\begin{conj}
    Let $\Omega$ be a bounded convex domain in $\mathbb{C}$, and consider a continuous function $h\colon\mathbb{C}\to\mathbb{R}$ that is constant outside of $\Omega$. Then for every $\theta\in (0,\pi)$ and positive $a$, there is a rectangle of radius $a$ centered in $\Bar{\Omega}$ with angle $\theta$ between its two diagonals such that $h$ takes the same value on its vertices.
\end{conj}

\section*{AI Declaration}
No AI tools were used to write this article, and all the ideas were developed by the author.

\bibliographystyle{plain}
\bibliography{ref}

\end{document}